\documentclass[preprint,1p,11pt]{elsarticle}
\usepackage{mathrsfs}
\usepackage{epsfig}
\usepackage{amsmath}
\usepackage{amssymb}
\usepackage{slashed}
\newcommand{\bea}{\begin{eqnarray}}
\newcommand{\eea}{\end{eqnarray}}
\newcommand{\bnn}{\begin{eqnarray*}}
\newcommand{\enn}{\end{eqnarray*}}
\newcommand{\be}{\begin{equation}}
\newcommand{\ee}{\end{equation}}

\journal{Advances in Applied Clifford Algebras}
\def\PACS{\par\leavevmode\hbox {\it PACS:\ }}%
\def\MSC{\par\leavevmode\hbox {\it MSC:\ }}%

\begin{document}
\begin{frontmatter}

\title{Conformal numbers}

\author{S. Ulrych}
\address{Wehrenbachhalde 35, CH-8053 Z\"urich, Switzerland}

\begin{abstract}
The conformal compactification is considered in a hierarchy of hypercomplex projective spaces with relevance in physics 
including Minkowski and Anti-de Sitter space.
The geometries are expressed in terms of bicomplex Vahlen matrices and further broken down into their structural components.
The relation between two subsequent projective spaces is displayed in terms of the complex unit and three additional
hypercomplex numbers. 
\end{abstract}



\end{frontmatter}

\section{Introduction}
\label{intro}
The complex numbers are central for the representation of physical processes in terms of mathematical models.
However, they cannot cover all aspects alone and generalizations are necessary. One of the possible generalizations
with potentially underestimated relevance in physics is the bicomplex number system \cite{Fut28,Sco34,Spa36,Ril53,Pri91}. 
Recent investigations in this area have been provided for example by \cite{Lun13,Str13,Che15,Kim16,Mur16}.
The number system is also known under the name of Segre numbers \cite{Seg92}.
The bicomplex numbers coincide with the
combination of hyperbolic and complex numbers formed by two commutative imaginary units, 
$i\equiv\sqrt{-1}$ and $j\equiv\sqrt{+1}$. The hyperbolic unit carries here the second complex number of
the bicomplex number system. More details on hyperbolic numbers have been provided beside many other authors by Yaglom \cite{Yag79}, Sobczyk \cite{Sob95}, Gal \cite{Gal02},
or in correspondence with the bicomplex numbers by Rochon and Shapiro \cite{Roc04}. 
The hypercomplex number systems are strongly connected to the theory of Clifford algebras and Lie groups, see Porteous \cite{Por95} or Ab{\l}amowicz et al. \cite{Abl96,Abl04}.
Algebraic properties of higher dimensional geometric spaces can be investigated in terms of hypercomplex matrix representations
of Clifford algebras. The generalizations can be applied also to functional calculus.
The properties of holomorphic functions of one complex variable can be extended to functions with values
in a Clifford algebra, consider here for example Brackx, Delanghe, and Sommen \cite{Bra82}. Further details can be found
also in G\"urlebeck et al. \cite{Gue08} or Colombo et al. \cite{Col04}.

The most prominent hypercomplex number system, representing the Clifford algebra $\mathbb{R}_{0,2}$, is given by the quaternions.
Girard shows in a short summary that quaternions provide spin representations 
of the most important equations in quantum physics and classical field theory \cite{Gir84}, see also Majern\'ik and Nagy \cite{Maj76}.
An extension is given by the spinor of conformal space, 
the twistor with its inherent null-plane geometry introduced by Penrose \cite{Pen67}.
Mathematically such an extension is described in terms of complex projective spaces \cite{Gol99}.
In the sixties of the last century conformal physics became popular in the context of the strong interaction. 
Conformal field theories play a central role within string
theories and the AdS/CFT (Anti-de Sitter space/conformal field theory) correspondence
of Maldacena \cite{Mal98}, Gubser et al. \cite{Gub98}, and Witten \cite{Wit98}.
In this area is situated also the higher spin holography \cite{Gab13,Ahn15,Vas15, Var15}, which could potentially benefit from the spin representations discussed in the following sections.
M\"obius geometries, conformal transformations, and their action on cycles
have been studied with focus on hypercomplex variables by Kisil \cite{Kis12}. 
Conformal extensions in relation with Clifford algebras and quaternionic analysis 
have been studied beside others by Sobczyk \cite{Sob12} and Frenkel and Libine \cite{Fre15}.

The mentioned publications provide the context for
a work about conformal relativity with hypercomplex variables, which has been published recently \cite{Ulr14}.
The motivation for a follow-up article was initially the consolidation of the method discussed in \cite{Ulr14}.
The objective was to provide easy to use mathematical tools within conformal physics in order to proceed
with the calculation of scattering amplitudes for comparison with experiment. From a conceptual point of view
it turned out that important geometries in physics, like the plane, Minkowski space, and the Anti-de Sitter space,
should be aligned and connected with a common set of rules.
Such connections between different geometries have been described in mathematical generality in terms of 
the conformal compactification with $2\times 2$ Clifford matrices,
a method which traces back to Vahlen \cite{Vah02} and Ahlfors \cite{Ahl85,Ahl86}.
Consider here also Porteous \cite{Por95} or Hertrich-Jeromin \cite{Her03}.
Maks \cite{Mak89} investigated explicitly the hierarchy
of M\"obius geometries, which will be used also in the following representation.
Each level in this hierarchy introduces an additional hypercomplex projective space based on the preceding geometry.  
On all levels M\"obius transformations can be applied to the considered hypercomplex variables.
The hypercomplex projective line and M\"obius spaces coincide throughout this hierarchy, which is in its complex restriction
only the case within two-dimensional geometries \cite{Sha97}. A similar generalization for quaternions is
of relevance in physics with respect to the concept of instantons, see Atiyah and Ward \cite{Ati77}.

The hierarchy of M\"obius geometries is introduced in the following sections based on the complexified null-plane numbers \cite{Ulr14}.
These numbers refer to the idempotents, which appear typically in the context of hyperbolic numbers.
The hierarchy of projective geometries begins with the two-dimensional plane as a non-trivial base manifold represented 
by complex numbers. This brings the method in \cite{Ulr14} 
closer to the twistor programme with its interpretation of space-time points as derived objects \cite{Pen77}. 
Furthermore, it reminds of the dimensional reduction of Minkowski space discussed by 't Hooft \cite{Hoo93}.
These investigations led to what is known today as the holographic principle, see also Susskind \cite{Sus95}.
It has been suggested to consider the holographic principle as a foundation of a quantum gravity theory,
in the same way as the equivalence principle is a foundation of general relativity \cite{Bou02}.
One may understand the hierarchy of M\"obius geometries also from the perspective of AdS/CFT.
General relativity on, e.g., AdS$_5$ is dual to the conformal field theory on its conformal boundary, which is given by Minkowski space.
The isometries within AdS$_5$ act as conformal transformations in Minkowski space \cite{Blu13}.

\section{Complex numbers}
\label{complex}
As mentioned in the introduction the two-dimensional Euclidean plane, represented by the complex numbers,
is the first non-trivial geometry to be considered. 
The intention of this section is also to provide simple examples for the notation used in this article.
In the context of Clifford algebras
a complex number corresponds to a Clifford paravector,
which is expanded in terms of the basis $e_\mu=(1,e_k)$. In this case there is only the single non-trivial basis element
\be
\label{single}
e_1=i\,.
\ee
As usual $i$ denotes the complex unit. For the sake of completeness it is worth to note that
the imaginary unit squares to the negative identity element and changes sign 
with respect to conjugation
\be
\label{squarecomplex}
i^2= -1\,,\quad \bar{i}=-i\,.
\ee
The complex unit can be considered as the single basis element of the Clifford algebra $\mathbb{R}_{0,1}$.

Two fundamental products can be defined for the complex numbers, the real and the complex product.
Consider here Andreescu and Andrica \cite{And14} especially with respect to the analogy of these products to the scalar and the cross product.
The real product is defined by its action on the basis elements $e_\mu=(1,e_k)$ of the paravector algebra
\be
\label{symprod}
e_\mu\cdot e_\nu=\frac{e_\mu\bar{e}_\nu+e_\nu\bar{e}_\mu}{2}=g_{\mu\nu}\,.
\ee
The metric tensor $g_{\mu\nu}$ has been introduced in this equation.
Inserting the definition of the paravector one finds the explicit result
\be
g_{\mu\nu}=\left(\begin{array}{cc}
1&0\\
0&1\\
\end{array}\right).
\ee
As expected the metric of the plane $\mathbb{R}^{2,0}$ is obtained.

The complex product is defined in analogy to the real product, but with a negative sign between the two contributions
\be
\label{sigma}
e_\mu\wedge e_\nu=\frac{e_\mu\bar{e}_\nu-e_\nu\bar{e}_\mu}{2}=\sigma_{\mu\nu}\,.
\ee
The anti-symmetric tensor $\sigma_{\mu\nu}$ has been introduced, which will be denoted in the following as spin tensor.
An explicit calculation based on the preceding definitions is leading to the following result
\be
\label{singlesigma}
\sigma_{\mu\nu}=\left(\begin{array}{cc}
0&-e_1\\
e_1&0\\
\end{array}\right).
\ee
The spin tensor is directly related to the spin angular momentum operator, where the spin tensor is just divided by a factor of two
\be
\label{spin}
s_{\mu\nu}=\frac{\sigma_{\mu\nu}}{2}\,.
\ee
In the context of Clifford algebras rotations can be defined with  the spin angular momentum operator in
a notation, which is still valid in higher dimensional spaces
\be
\label{rota}
r=\exp{\left(\frac{1}{2}s_{\mu\nu}\omega^{\mu\nu}\right)},\quad
\omega^{\mu\nu}=-\omega^{\nu\mu}.
\ee
The rotation is acting on a complex number by $rzr^\dagger$ or equivalently by $rz\hat{r}^{-1}$, where
$r^\dagger$ indicates reversion and $\hat{r}$ the main involution \cite{Por95,Gue08}.

\section{Complex null-plane}
\label{null plane}
The complex null-plane numbers have appeared in \cite{Ulr14} as a key structure
in a representation of M\"obius geometries based on
hyperbolic and complex units. 
In a real form null-plane numbers were considered before for example by Zhong \cite{Zho85} and Hucks \cite{Huc93}
in the context of the hyperbolic number system.  The complexified null-plane numbers can be defined by the following equations
\be
\label{nullcalc}
oo=i o\,,\quad \bar{o}\bar{o}=-i\bar{o}\,,\quad o\bar{o}=0\,.
\ee
The bar symbol denotes conjugation. One may express the null-plane numbers in an alternative 
notation with the complex unit and a second hypercomplex structure
\be
\label{derive}
i = o-\bar{o}\,,\quad j=o+\bar{o}\,.
\ee
It should be noted that the hypercomplex unit $j$ 
corresponds to $ij$ in the notation of \cite{Ulr14}\,.
The square of $j$  can be calculated with the preceding definitions 
\be
\label{squares}
j^2= -1\,,\quad \bar{j}=j\,.
\ee
The two imaginary units $i$ and $j$ show a different behaviour with respect to conjugation, which can be derived from Eq.~(\ref{derive})
with the rule $\bar{\bar{o}}=o$.
The imaginary unit $j$ is invariant with respect to conjugation,
whereas the complex unit $i$ changes sign. The combination of these units results in an algebra, which is among the six representations of
bicomplex numbers investigated by Alpay et al. \cite{Alp14}.

\section{Conformal compactification}
\label{compact}
The complex space is by its nature unlimited. However, with a conformal compactification it is possible to enclose the
unlimited geometry by adding infinity. This leads to the projective space $\mathbb{P}^{1}_\mathbb{C}=\mathbb{C}\cup\{\infty\}$.
In order to perform the compactification the algebra introduced in the preceding section has to be
multiplied by additional $2\times 2$ matrix structures, see for example Obolashvili \cite{Obo98}.
Therefore two additional units are introduced, which are based on explicit matrix representations
\be
\imath=\left(\begin{array}{cc}
0&1\\
-1&0
\end{array}\right),\quad \jmath=\left(\begin{array}{cc}
0&1\\
1&0
\end{array}\right).
\ee
The matrix $\imath$ can be seen as the counterpart of $i$ as it changes sign
under conjugation, which is represented by transposition of the matrix.
In contrast, $\jmath$ remains invariant with respect to conjugation, but squares to the identity element.
Multiplication of these matrices results in
\be
\imath\jmath=\left(\begin{array}{cc}
1&0\\
0&-1
\end{array}\right)=-\jmath\hspace{0.03cm}\imath\,.
\ee
The three matrices correspond to the Lie algebra of $SL(2,\mathbb{R})$.

With these matrices and the complex null-plane numbers the conformal compactification can be represented more 
compact and generalized compared to \cite{Ulr14}.
The base geometry is considered to have an even number of $2m=n$ dimensions.
The $n-1$ basis elements of the Clifford algebra $\mathbb{R}_{m-1,m}$ are transformed to the basis elements of the projective space by
\be
\label{basspace}
e_{k}=\imath\jmath e_k\,,\quad k=1,\dots,n-1\,.
\ee
On the right hand side of the equation are the basis elements of the source geometry.
The two additional basis elements of the conformal algebra are given by
\be
\label{basrest}
e_{n}=i \jmath
\,,\quad
e_{n+1}=\imath j\,.
\ee
The resulting basis elements generate the Clifford algebra $\mathbb{R}_{m,m+1}$. 
The argument to enclose the original space by adding infinity applies to arbitrary even dimensional spaces in this hierarchy.
Thus there is an infinite series of conformal compactifications. One may understand the infinite series of M\"obius geometries
as representation spaces of a dimension independent projective scheme.

\section{Minkowski space}
\label{Minkowski}
The method discussed in the previous section can now be applied explicitly to
the complex numbers. With Eqs.~(\ref{basspace}) and (\ref{basrest}) one can derive the following three basis elements, which can 
be used to introduce a paravector model $e_\mu=(1,e_k)$ of Minkowski space
\be
\label{Pauli}
e_1=\left(\begin{array}{cc}
i&0\\
0&-i
\end{array}\right),\quad
e_2=\left(\begin{array}{cc}
0&i\\
i&0
\end{array}\right),\quad
e_3=\left(\begin{array}{cc}
0&j\\
-j&0
\end{array}\right).
\ee
The algebra corresponds to the Clifford algebra $\mathbb{R}_{1,2}$
and is thus isomorphic to the Pauli algebra \cite{Por95}.
The Clifford algebra $\mathbb{R}_{1,2}$ as introduced above will replace the algebra $\mathbb{R}_{3,0}$ in \cite{Ulr14}.

The metric tensor can be calculated with Eq.~(\ref{symprod}) using
the basis elements introduced above
\be
\label{relmetric}
g_{\mu\nu}=\left(\begin{array}{cccc}
1&0&0&0\\
0&1&0&0\\
0&0&1&0\\
0&0&0&-1\\
\end{array}\right).
\ee
The metric convention of Minkowski space $\mathbb{R}^{3,1}$ has been reversed compared to \cite{Ulr14}.
The spin tensor $\sigma_{\mu\nu}$ can be calculated with Eq.~(\ref{sigma})
\be
\label{ordinaryspin}
\sigma_{\mu\nu}=\left(\begin{array}{cccc}
0&-e_1&-e_2&-e_3\\
e_1&0&-j e_3&-j e_2\\
e_2&j e_3&0&j e_1\\
e_3&j e_2&-j e_1&0\\
\end{array}\right).
\ee
The spin angular momentum operator is given again by Eq.~(\ref{spin}).
With the multiplication rules of the hypercomplex variables the commutation relation
of the relativistic spin angular momentum can be calculated
\be
\label{comm}
[s_{\mu\nu},s_{\rho\sigma}]=g_{\mu\sigma}s_{\nu\rho}-g_{\mu\rho}s_{\nu\sigma}+g_{\nu\sigma}s_{\mu\rho}
-g_{\nu\rho}s_{\mu\sigma}\,.
\ee
Thus the spin matrices provide a representation of the Lorentz group $SO(3,1,\mathbb{R})$. 
In comparison to \cite{Ulr14} the time coordinate has to switch to $e_3$.
In this sense pure rotations are represented within the paravector model $e_\mu=(1,e_1,e_2)$.
One can see in the matrix representations of the spin tensor, that 
the rotations are free of the imaginary unit $j$. Boosts and time thus come in relation with $j$.

This becomes more obvious if one breaks up the above $SO(3,1,\mathbb{R})$ spin
representation with Eqs.~(\ref{basspace}) and (\ref{basrest}) in terms of the single basis element of the complex numbers, $e_1=i$.
This leads to the spin tensor
\be
\sigma_{\mu\nu}=\left(\begin{array}{cccc}
0&-\imath\jmath e_1&-i\jmath&-\imath j\\
\imath\jmath e_1&0&-\imath i e_1&-\jmath j e_1\\
i\jmath&\imath i e_1&0&\imath\jmath ij\\
\imath j&\jmath je_1&-\imath\jmath ij&0\\
\end{array}\right).
\ee
This representation is equivalent to Eq.~(\ref{ordinaryspin}).
The time coordinate $t$ is attached to $e_3$ as mentioned before.
The space dimensions $(x,y,z)$ can be attributed in arbitrary rotated form to $e_\mu=(1,e_1,e_2)$. The base
algebra $(1,e_1)$ is still included. It can be interpreted as a geometric polarization plane
and is potentially related to interacting forces and masses.
In order to indicate a possible relation between geometry and particle masses one can introduce the following equation, 
which is inspired by the Regge trajectories in the sense that the angular momentum is related to squared masses \cite{Gui91}
\be
\label{proton}
4\pi\exp{(4\pi)}=\left(\frac{m_p}{m_e}\right)^2\,.
\ee
The experimental proton to electron mass ratio is calculated with a deviation of $3.4\%$. The question arises whether this equation can be 
derived from geometric properties of manifolds, which describe single protons and electrons \cite{Ati12}. 

\section{Anti-de Sitter space}
\label{confalg}
One can apply the conformal compactification to Minkowski space
and reaches the ambient space $\mathbb{R}^{4,2}$, which includes the Anti-de Sitter space AdS$_5$.
Equations~(\ref{basspace}) and (\ref{basrest}) are applied to the basis elements of $\mathbb{R}_{1,2}$ representing
the base manifold.
The resulting new basis elements of the Clifford algebra $\mathbb{R}_{2,3}$ are used to set up the paravector model $e_\mu=(1,e_k)$. 
The metric tensor is calculated with Eq.~(\ref{symprod})
\be
g_{\mu\nu}=\left(\begin{array}{cccccc}
1&0&0&0&0&0\\
0&1&0&0&0&0\\
0&0&1&0&0&0\\
0&0&0&-1&0&0\\
0&0&0&0&1&0\\
0&0&0&0&0&-1\\
\end{array}\right).
\ee
The spin tensor is computed with Eq.~(\ref{sigma}). The result
can be displayed in terms of the three basis elements of the Clifford algebra $\mathbb{R}_{1,2}$ and the four
hypercomplex units, which have been introduced in the preceding sections
\be
\label{confspin}
\sigma_{\mu\nu}=\left(\begin{array}{cccccc}
0&-\imath\jmath e_1&-\imath\jmath e_2&-\imath\jmath e_3&-i\jmath&-\imath j\\
\imath\jmath e_1&0&-j e_3&-j e_2&-\imath ie_1&-\jmath je_1\\
\imath\jmath e_2&j e_3&0&j e_1&-\imath ie_2&-\jmath je_2\\
\imath\jmath e_3&j e_2&-j e_1&0&-\imath ie_3&-\jmath je_3\\
i\jmath&\imath ie_1&\imath ie_2&\imath ie_3&0&\imath\jmath ij\\
\imath j&\jmath je_1&\jmath je_2&\jmath je_3&-\imath\jmath ij&0\\
\end{array}\right).
\ee
The spin angular momentum operator is defined again by Eq.~(\ref{spin})
and satisfies the commutation relations of Eq.~(\ref{comm}).
The first column in the spin tensor includes the basis elements of the Clifford algebra $\mathbb{R}_{2,3}$,
which will replace $\mathbb{R}_{4,1}$ in \cite{Ulr14}.
The algebra is equivalent to the Dirac algebra \cite{Por95}.

\section{Conformal spin in the base space}
\label{minkconf}
The elements of the spin algebra of an $n+2$ dimensional ambient space, 
can be represented within the $n$ dimensional base manifold.
For example the spin angular momentum defined by Eqs.~(\ref{spin}) and (\ref{confspin})
can be restricted to indices $\mu=0,\dots,n-1$, in this case with $n=4$. The
reduced spin angular momentum operator then still satisfies Eq.~(\ref{comm}) with the corresponding metric tensor.
The remaining operators can be reorganized using the definitions of Kastrup \cite{Kas62},
which result in spin representations of the conformal group
\bea
\label{coprot}
p_\mu&=&-s_{\mu n}-s_{\mu n+1}\,,\nonumber\\
q_\mu&=&s_{\mu n}-s_{\mu n+1}\,,\nonumber\\
d&=&s_{nn+1}\,.
\eea
Here $p_\mu$ labels the spin representation of the momentum operator.
The notation $q_\mu$ is used for the spin representation of
the special conformal transformations and $d$ for the scale transformation. 

Based on these definitions explicit expressions for the spin representation of these operators
can be calculated 
\bea
\label{cop}
2p_0&=&\imath j+i\jmath\,,\nonumber\\
2q_0&=&\imath j-i\jmath\,,\nonumber\\
2p_k&=&e_k(\jmath j +\imath i)\,,\nonumber\\
2q_k&=&e_k(\jmath j -\imath i)\,,\nonumber\\
2d&=&\imath\jmath ij\,.
\eea
The basis elements $e_k$ refer to the base manifold, this means to the right hand side of Eq.~(\ref{basspace}).
Beside the already mentioned commutation relation of the spin angular momentum 
one finds the following commutation relations
\bea
\label{confcomm}
\left[s_{\mu\nu},p_\sigma\right]&=&g_{\nu \sigma}p_{\mu}-g_{\mu \sigma}p_{\nu}\,,\nonumber\\
\left[s_{\mu\nu},q_\sigma\right]&=&g_{\nu \sigma}q_{\mu}-g_{\mu \sigma}q_{\nu}\,,\nonumber\\
\left[d,p_\mu\right]&=&-p_\mu\,,\nonumber\\
\left[d,q_\mu\right]&=&q_\mu\,,\nonumber\\
\left[q_\mu,p_\nu\right]&=&2(g_{\mu\nu}d+s_{\mu\nu})\,.
\eea
All other commutators vanish. Thus the spin operators form a representation of the conformal group.
Due to the relation between conformal transformations and M\"obius geometries, it should be noted that
the M\"obius space is identified with the homogeneous space defined by the Lie algebras
\be
\mathfrak{g}=\{s_{\mu\nu},p_\mu,q_\mu,d\}\,,\quad\mathfrak{h}=\{s_{\mu\nu},q_\mu,d\}\,.
\ee
More detailed information about M\"obius geometries can be found in the textbook of Sharpe \cite{Sha97}.

With respect to physics one finds that the
Minkowski space is situated within the sequence of projective geometries.
Therefore Eq.~(\ref{coprot}) points to the dimensional reduction considered by 't Hooft \cite{Hoo93} and furthermore to the AdS/CFT correspondence.
The isometries on a sphere in Minkowski space $\mathbb{R}^{3,1}$ are dual to 
conformal transformations in the Euclidean planar limit \cite{Hoo73}. 
The conformal transformations refer to the symmetries of the Laplace equation in $\mathbb{R}^{2,0}$, 
which can be chosen to define the conformal field theory
\be
\triangle\psi=0\,.
\ee
The solutions give rise to geometric fields, which can be transformed with the above symmetry operators represented in the corresponding function space.

\section{Summary}
A system of hypercomplex units has been defined for the representation
of relativistic physics in terms of paravector models. The representation is used to attach
a light cone to a given base manifold by virtue of the conformal compactification.
The method is applied to the complex numbers representing the initial non-trivial base space,
which results in a hypercomplex representation of Minkowski space. 

The system of hypercomplex units is used furthermore to create hypercomplex spin representations within
M\"obius geometries.
This is leading to a finite dimensional representation of the conformal symmetry operators of the
two dimensional Laplace equation. Higher dimensional geometries
of physical relevance can be derived around the base manifold, which can be accessed through
holomorphic functions of ordinary complex numbers.

\end{document}